\documentclass[12pt]{amsart}

\title{A going-up theorem}
\author{Ying Zong}
\address{Department of Mathematics\\University of Toronto} \email{zongying@math.utoronto.ca}

\date{} 

\begin{document}

\maketitle


Let $p$ be a prime number and let $\Lambda=\mathbf{Z}/p\mathbf{Z}$. For every scheme $X$ and for $*\in\{+, -, b, \emptyset\}$, let $D^*(X, \Lambda)$ denote the derived category with the boundedness condition $*$ of the abelian category of sheaves of $\Lambda$-modules on $(\mathrm{Sch}/X)$ for the \'{e}tale topology. 
\smallskip

The following result, which we call \emph{going-up}, is the consequence of a reading of Zariski's paper \emph{On the irregularity of cyclic multiple planes} and has only to do with certain part of his extremely attractive algebro-topological arguments.

\smallskip

\smallskip

{\bf Theorem.} --- \emph{Let $X$ be a scheme, let $D$ be a closed sub-scheme of $X$ and let $U$ be the complement of $D$ in $X$. Assume that $U$ is locally connected. Let $f: Y\to X$ be an $X$-scheme which is entier surjective over $X$, radical over $D$ and finite \'{e}tale Galois over $U$ with Galois group a $p$-group. Let
\[T: D^*(X, \Lambda)\to A\] be a cohomological functor with values in an abelian category $A$. Assume that $K$ is an object of $D^*(X, \Lambda)$ such that $T(K[1])=0$ and that the natural morphism
\[T(K)\to T(i_*i^*K)\] is an epimorphism, where $i: D\hookrightarrow X$ denotes the canonical closed immersion.}
\smallskip

\emph{Then $T(f_*f^*K[1])=0$.}  

\begin{proof} Let $j: U\hookrightarrow X$ be the open immersion. With each object $L$ of $D^*(X, \Lambda)$ one associates a localization triangle (SGA 4 XVII 5.1.16)
\[j_{!}j^*L\to L\to i_*i^*L,\] from which one deduces the $T$-cohomology sequence
\[\cdots\to T^1(j_{!}j^*L)\to T^1(L)\to T^1(i_*i^*L)\stackrel{\delta_L}{\longrightarrow} T^2(j_{!}j^*L)\to\cdots\] where $T^n(N)$ stands for $T(N[n])$ for every integer $n$ and every object $N$ of $D^*(X, \Lambda)$. 
\smallskip

In particular, when considering $L=f_*f^*K$, one can infer that 
$T^1(L)=0$ provided that 
\[T^1(j_{!}j^*L)=0\] and that
\[\delta_L: T^1(i_*i^*L)\to T^2(j_{!}j^*L)\] is a monomorphism.
\smallskip

\smallskip

--- \emph{For $L=f_*f^*K$, one has $T^1(j_{!}j^*L)=0$ }:
\smallskip

\smallskip

The restriction of $f$ to $f^{-1}(U)$, $f^{-1}(U)\to U$, is by assumption finite \'{e}tale Galois with Galois group a $p$-group $G$. The $\Lambda$-module $j^*f_*\Lambda$ on $U$ is therefore locally constant constructible and its corresponding monodromy representation at each geometric point $u$ of $U$,
\[\rho_u: \pi_1(U, u)\to \mathrm{GL}_{\Lambda}((j^*f_*\Lambda)_u),\] has as monodromy a subquotient of $G$. Thus this monodromy $\mathrm{Im}\ \rho_u$ is a unipotent subgroup of $\mathrm{GL}_{\Lambda}((j^*f_*\Lambda)_u)$, and so, if one sets
\[F^1:=\mathrm{Coker}\ (\Lambda\to j^*f_*\Lambda),\] the $\Lambda$-module $F^1$ admits a finite filtration by sub-$\Lambda$-modules $F^*$ such that each successive quotient $gr^*_F$ is isomorphic to $\Lambda$. One applies now the triangulated functor 
\[-\mapsto j^*K\stackrel{\mathbf{L}}{\otimes}_{\Lambda} -\] to the exact sequences
\[0\to \Lambda\to j^*f_*\Lambda\to F^1\to 0\]
\[0\to F^{*+1}\to F^*\to \Lambda\to 0.\]

This gives rises to the triangles 
\[j^*K\to j^*K\stackrel{\mathbf{L}}{\otimes}_{\Lambda}j^*f_*\Lambda\to j^*K\stackrel{\mathbf{L}}{\otimes}_{\Lambda}F^1\]
\[j^*K\stackrel{\mathbf{L}}{\otimes}_{\Lambda}F^{*+1}\to j^*K\stackrel{\mathbf{L}}{\otimes}_{\Lambda}F^*\to j^*K\] to which one further applies the cohomological functor $Tj_{!}$ to obtain the exact sequences 
\[T^1(j_{!}j^*K)\to T^1(j_{!}(j^*K\stackrel{\mathbf{L}}{\otimes}_{\Lambda}j^*f_*\Lambda))\to T^1(j_{!}(j^*K\stackrel{\mathbf{L}}{\otimes}_{\Lambda}F^1))\]
\[T^1(j_{!}(j^*K\stackrel{\mathbf{L}}{\otimes}_{\Lambda}F^{*+1}))\to T^1(j_{!}(j^*K\stackrel{\mathbf{L}}{\otimes}_{\Lambda}F^*))\to T^1(j_{!}j^*K).\]

\smallskip

Observe that 
\smallskip

\smallskip

--- \emph{One has $T^1(j_{!}j^*K)=0$ }:
\smallskip

\smallskip

Indeed, the exact sequence
\[T(K)\to T(i_*i^*K)\to T^1(j_{!}j^*K)\to T^1(K)\] shows that this follows from the hypotheses that  
\[T^1(K)=0\] and that
\[T(K)\to T(i_*i^*K)\] is an epimorphism.

\smallskip

\smallskip

Hence, by descending induction on $n\geq 1$, one obtains that
\[T^1(j_{!}(j^*K\stackrel{\mathbf{L}}{\otimes}_{\Lambda}j^*f_*\Lambda))=0\]
\[T^1(j_{!}(j^*K\stackrel{\mathbf{L}}{\otimes}_{\Lambda}F^n))=0\] for all $n\geq 1$.
\smallskip

Finally, using the projection formula (SGA 4 XVII 5.2.9)
\[f_*f^*K=K\stackrel{\mathbf{L}}{\otimes}_{\Lambda}f_*\Lambda\] one concludes that
\[T^1(j_{!}j^*f_*f^*K)=T^1(j_{!}j^*(K\stackrel{\mathbf{L}}{\otimes}_{\Lambda}f_*\Lambda))=T^1(j_{!}(j^*K\stackrel{\mathbf{L}}{\otimes}_{\Lambda}j^*f_{*}\Lambda))=0.\]

\smallskip

\smallskip

--- \emph{For $L=f_*f^*K$, the coboundary morphism
\[\delta_L: T^1(i_*i^*L)\to T^2(j_{!}j^*L)\] is a monomorphism }:

\smallskip

\smallskip

Let $i': f^{-1}(D)\hookrightarrow Y$ be the canonical closed immersion and let $f': f^{-1}(D)\to D$ be the restriction of $f$ to $f^{-1}(D)$. Notice that (SGA 4 VIII 5.5) one has  
\[i^*f_{*}=f'_{*}i^{'*}\] and hence
\[T^1(i_*i^*f_*f^*K)=T^1(i_*f'_{*}i^{'*}f^*K)=T^1(i_*f'_{*}f^{'*}i^*K).\] 

Notice next that the adjunction morphism
\[i^*K\to f'_{*}f^{'*}i^*K\] induced by the pair of adjoint functors $f'_{*}, f^{'*}$ is an isomorphism, for $f'$ is entier radical surjective (SGA 4 VIII 1.1).
\smallskip

So one has the identity
\[T^1(i_*i^*f_*f^*K)=T^1(i_*i^*K),\] by which the coboundary morphism
\[\delta_L: T^1(i_*i^*L)\to T^2(j_{!}j^*L),\] for $L=f_*f^*K$, can be identified with the composition
\[T^1(i_*i^*K)\stackrel{\delta_K}{\longrightarrow}T^2(j_{!}j^*K)\to T^2(j_{!}j^*f_*f^*K).\] 

It suffices to show that both
\[\delta_K: T^1(i_*i^*K)\to T^2(j_{!}j^*K)\] and
\[T^2(j_{!}j^*K)\to T^2(j_{!}j^*f_{*}f^*K)\] are monomorphisms.

\smallskip

\smallskip

--- \emph{The coboundary
\[\delta_K: T^1(i_{*}i^*K)\to T^2(j_{!}j^*K)\] is a monomorphism }:
\smallskip

\smallskip

For, $T^1(K)=0$ and one has the exact sequence
\[T^1(K)\to T^1(i_*i^*K)\to T^2(j_{!}j^*K).\]

\smallskip

--- \emph{The morphism
\[T^2(j_{!}j^*K)\to T^2(j_{!}j^*f_*f^*K)\] induced by the pair of adjoint functors $f_*, f^*$ is a monomorphism }:

\smallskip

\smallskip

The triangle
\[j^*K\to j^*K\stackrel{\mathbf{L}}{\otimes}_{\Lambda}j^*f_*\Lambda\to j^*K\stackrel{\mathbf{L}}{\otimes}_{\Lambda}F^1\] gives rises to the following cohomology sequence
\[T^1(j_{!}(j^*K\stackrel{\mathbf{L}}{\otimes}_{\Lambda}F^1))\to T^2(j_{!}j^*K)\to T^2(j_{!}(j^*K\stackrel{\mathbf{L}}{\otimes}_{\Lambda}j^*f_{*}\Lambda)).\]

It remains only to recall that
\[j_{!}(j^*K\stackrel{\mathbf{L}}{\otimes}_{\Lambda}j^*f_*\Lambda)=j_{!}j^*f_{*}f^*K\] and that
\[T^1(j_{!}(j^*K\stackrel{\mathbf{L}}{\otimes}_{\Lambda}F^1))=0.\]

\end{proof}

\smallskip

{\bf Remark.} --- The following analogue is proven in the same way :
\smallskip

\emph{Let $X$ be a topological space, let $D$ be a closed subspace of $X$ and let $U$ be the complement of $D$ in $X$. Assume that $U$ is locally path simply connected. Let $f: Y\to X$ be an $X$-space which is proper surjective with finite fibers over $X$ and which is a bijection above $D$ and is a Galois covering with Galois group a $p$-group above $U$. Let 
\[T: D^*(X, \mathbf{Z}/p\mathbf{Z})\to A\] be a cohomological functor with values in an abelian category $A$. Assume that $K$ is an object of $D^*(X, \mathbf{Z}/p\mathbf{Z})$ such that $T(K[1])=0$ and that the natural morphism
\[T(K)\to T(i_*i^*K)\] is an epimorphism, where $i: D\hookrightarrow X$ denotes the canonical inclusion.}
\smallskip

\emph{Then $T(f_*f^*K[1])=0$.}

\smallskip

\smallskip

\smallskip

\smallskip

Consider $T=H^0\mathrm{R}\Gamma: D^+(X, \Lambda)\to (Ab)$, where $\Lambda=\mathbf{Z}/p\mathbf{Z}$.
\smallskip

Take $K=\Lambda$. Assume $H^1(X, \Lambda)=0$ and that $H^0(X, \Lambda)\to H^0(D, \Lambda)$ is surjective. Then one obtains that $H^1(Y, \Lambda)=0$.
\smallskip

Or, take $K=j_{!}\Lambda[n-1]$ for an integer $n$. As $i^*K=0$, one finds that $H^n(Y, j'_{!}\Lambda)=0$ when $H^n(X, j_{!}\Lambda)=0$. Here $j': f^{-1}(U)\hookrightarrow Y$ is the inclusion. 
\smallskip

\smallskip

\smallskip

\smallskip


\bibliographystyle{amsplain}

\begin{thebibliography}{1}

\bibitem{}
O.~Zariski.
\newblock On the irregularity of cyclic multiple planes.
\newblock \emph{Annals of Mathematics}, 32, 1931.



\end{thebibliography}

\end{document}